# Counting dimensions of $L$-harmonic functions

By Peter Li and Jiaping Wang*

## 0. Introduction

In this article, we will consider second order uniformly elliptic operators of divergence form defined on $\mathbb{R}^n$ with measurable coefficients. Mainly, we will give estimates on the dimension of space of solutions that grow at most polynomially of degree $d$. More precisely, in terms of a rectangular coordinate system $\{x_1, \ldots, x_n\}$, a second order uniformly elliptic operator of divergence form, $L$, acting on a function $f \in H^1_{\text{loc}}(\mathbb{R}^n)$ is given by

$$Lf = \sum_{i,j=1}^{n} \frac{\partial}{\partial x_i}\left(a^{ij}(x)\frac{\partial f}{\partial x_j}\right),$$

where $(a^{ij}(x))$ is an $(n \times n)$-symmetric matrix satisfying the ellipticity bounds

(0.1) $$\lambda I \leq (a^{ij}) \leq \Lambda I$$

for some constants $0 < \lambda \leq \Lambda < \infty$. Other than the ellipticity bounds, we only assume that the coefficients $(a_{ij})$ are merely measurable functions.

A function $f \in H^1_{\text{loc}}(\mathbb{R}^n)$ is said to be $L$-harmonic if $Lf = 0$. Since $(a_{ij})$ are only measurable, this equation is being interpreted in the weak sense. The celebrated De Giorgi-Moser-Nash theory asserts that any weak solution $f$ must be $C^\alpha$ for some $0 < \alpha < 1$. A global version of this theory implies that there exists $0 < \alpha < 1$ such that if $\rho(x)$ denotes the Euclidean distance from $x$ to the origin, then any $L$-harmonic function $f$ satisfying the growth condition

$$|f(x)| = O(\rho^\alpha(x)),$$

*The research of the first author was partially supported by NSF grant #DMS-9626310. The research of the second author was partially supported by NSF grant #DMS-9704482.

Part of this work was done when the first author was visiting the Institute of Mathematical Sciences at the Chinese University of Hong Kong. He would like to thank the members of the Department of Mathematics for their hospitality during his visit.

as $x \to \infty$, must be identically constant. This statement can be thought of as a $C^\alpha$ regularity result at infinity. In general, for each real number $d \geq 0$, we can define

$$\mathcal{H}_d(L) = \{f \in H^1_{\mathrm{loc}}(\mathbb{R}^n) \,|\, Lf = 0 \text{ and } |f(x)| = O(\rho^d(x))\}$$

to be the space of all polynomial growth $L$-harmonic functions of degree at most $d$. Let us denote the dimension of $\mathcal{H}_d(L)$ by $h_d(L)$. The above growth result can be restated as

$$h_d(L) = 1$$

for all $0 \leq d \leq \alpha < 1$. The main purpose of this article is to give an upper bound estimate on $h_d(L)$.

Before stating our result precisely, let us give a brief outline on the history of this problem. After more than 20 years since the establishment of the De Giorgi-Moser-Nash theory, Avellaneda and Lin [A-Ln] studied the solution space of a uniformly elliptic operator of divergence form. They assumed that the coefficients are Lipschitz continuous and periodic. Using the theory of homogenization, they showed that every polynomial growth $L$-harmonic function is necessarily a polynomial of periodic coefficients. Moreover, there is a linear isomorphism between the harmonic polynomials and the polynomial growth $L$-harmonic functions on $\mathbb{R}^n$. In particular, this implies that $h_d(L)$ are finite for all $d$. Later, Moser and Struwe [M-S] using a different argument proved the same result without the Lipschitz continuity assumption on the coefficients (see also [Ln] for this last point). Their argument also works for a certain class of nonlinear operators. Recently, based on the concept of $G$-convergence, Lin [Ln] considered the case when $L$ is an asymptotically conic operator. In this case, his results imply that $h_d(L)$ are finite and have an upper bound depending on $n$, $d$, $\lambda$ and $\Lambda$. Zhang [Z] followed this direction by assuming that $L$ converges to a possibly $\mathbb{S}^1$ family of conic operators and established some similar results. Last year, the first author [L1] proved that for a general uniformly elliptic operator of divergence (and nondivergence form) with measurable coefficients, there exists a constant $C > 0$ depending only on $n$, $\lambda$ and $\Lambda$, such that,

$$(0.2) \qquad h_d(L) \leq C\, d^{n-1}$$

for all $d \geq 1$.

Note that when $L$ is the Euclidean Laplacian on $\mathbb{R}^n$, one can compute $h_d(\Delta)$ directly because $\mathcal{H}_d(\Delta)$ is spanned by the homogeneous harmonic poly-nomials of degree at most $d$. In fact, it is known that

$$h_d(\Delta) = \binom{n+d-1}{d} + \binom{n+d-2}{d-1}$$
$$= \frac{(2d+n-1)}{(n-1)!}(d+n-2)(d+n-3)\ldots(d+1)$$
$$= \frac{2}{(n-1)!}\left(d^{n-1} + O(d^{n-2})\right).$$

This example indicates that the estimate (0.2) is sharp in terms of the order of $d$. This computation also raises the question of whether it is possible to obtain an estimate for $h_d(L)$ which will reflect the sharp constant. One may ask if there exists a constant $C(n, \Lambda\lambda^{-1}) > 0$ depending only on the prescribed quantities and has the property that $C(n, \Lambda\lambda^{-1}) \to 1$ as $\Lambda\lambda^{-1} \to 1$, such that,

$$h_d(L) \leq C(n, \Lambda\lambda^{-1}) \frac{2}{(n-1)!} d^{n-1}.$$

In particular, such an estimate will be sharp as realized by $L = \Delta$.

This article is an attempt to prove the validity of this estimate. While we can not establish the estimate in the precise form, we do prove a sharp upper bound for the truncated sum of $h_d(L)$. In fact, our estimate only depends on the ratio of the ellipticity constants at infinity. More precisely, let us denote $\lambda_r$ and $\Lambda_r$ to be the ellipticity bounds satisfying

$$\lambda_r I \leq (a^{ij}(x)) \leq \Lambda_r I$$

for all $x \in \mathbb{R}^n \setminus B(r)$ with $B(r) = \{x \,|\, \rho(x) < r\}$ being the Euclidean ball of radius $r$. Obviously, $\lambda_r$ and $\Lambda_r$ are montonicially increasing and decreasing functions of $r$, respectively. Also, they satisfy

$$0 < \lambda \leq \lambda_r \leq \Lambda_r \leq \Lambda < \infty.$$

Let us define their limits at infinity by

$$\lambda_\infty = \lim_{r\to\infty} \lambda_r,$$

and

$$\Lambda_\infty = \lim_{r\to\infty} \Lambda_r.$$

The main result of this article is to prove that

$$\sum_{i-1}^{d} h_i(L) \leq \left(\frac{\Lambda_\infty}{\lambda_\infty}\right)^{n-1} \frac{2}{n!}(d^n + O(d^{n-1})).$$

This upper bound is sharp for $L = \Delta$ by simply summing up the formula for $h_d(\Delta)$.

We would like to point out that a similar type of estimate was proved by the authors [L-W2] for the Laplacian on a complete manifold with nonnegative sectional curvature. In both the nonnegatively curved manifold case and our present case, our methods rely on restricting the operator on codimension-1 submanifolds. A major difficulty for measurable coefficients is that the coefficients might not be measurable when restricted on a codimension-1 submanifold. As a result, extra work is needed to overcome this problem. In a parallel manner, the lack of regularity issue occurs on the cut-locus for the manifold case and extra but different effort was required to overcome that hurdle.

It is also worth mentioning that tremendous progress has been made in the direction of the establishment of Liouville properties and in counting the dimensions of polynomial growth harmonic functions on a manifold. We will refer the interested reader to the paper [C-M] by Colding and Minicozzi and a recent survey article by the first author [L2], and some subsequent articles [S-T-W], [L-W1].

## 1. $\phi$-harmonic functions on a complete manifold

Our strategy is to translate the problem of studying $L$-harmonic functions on $\mathbb{R}^n$ to that of studying $\phi$-harmonic functions on a manifold $(\mathbb{R}^n, g)$. Let $(M^n, g)$ be an $n$-dimensional Riemannian manifold. Let $\phi$ be a positive measurable function on $M$. A function $f$ is said to be $\phi$-harmonic if it satisfies

$$\mathrm{div}\,(\phi\,\nabla f) = 0,$$

where $\nabla$ and div are the gradient and the divergence with respect to the Riemannian metric $g$. In terms of local coordinates $\{x_1, \ldots, x_n\}$, the Riemannian metric can be expressed as $g = \sum_{i,j=1}^n g_{ij}\, dx_i\, dx_j$. Using this, we can write

$$(1.1) \qquad \mathrm{div}\,(\phi\,\nabla f) = \sum_{i,j=1}^n \frac{1}{\sqrt{G}} \frac{\partial}{\partial x_i} \left( \phi \sqrt{G} g^{ij} \frac{\partial f}{\partial x_j} \right),$$

where $G = \det(g_{ij})$ and $(g^{ij}) = (g_{ij})^{-1}$.

For any nonnegative real number $d$, we denote

$$\mathcal{H}_d(g, \phi) = \{f \mid \mathrm{div}\,(\phi\,\nabla f) = 0,\ |f(x)| = O(r^d(x))\}$$

to be the space of $\phi$-harmonic functions with at most polynomial growth of degree $d$ with respect to the geodesic distance $r(x)$ to some fixed point. Let us denote the dimension of $\mathcal{H}_d(g, \phi)$ by $h_d(g, \phi) = \dim \mathcal{H}_d(g, \phi)$.

To set up the $L$-harmonic problem as a $\phi$-harmonic problem, we recall that $\rho(x)$ is the Euclidean distance function to the origin. Let us introduce the Riemannian metric $g = \sum_{i,j=1}^n g_{ij}\, dx_i\, dx_j$ on $\mathbb{R}^n$ by

$$(1.2) \qquad (g^{ij}) = (g_{ij})^{-1} = w(x)\,(a^{ij}),$$

where

(1.3) $$w(x) = \sum_{i,j=1}^{n} a^{ij} \frac{\partial \rho}{\partial x_i} \frac{\partial \rho}{\partial x_j}.$$

If we define

(1.4) $$\phi_g(x) = w^{\frac{n-2}{2}}(x) \sqrt{\det(a^{ij})},$$

then one easily verifies by (1.1) that $f$ is $L$-harmonic if and only if it is $\phi_g$-harmonic on $M = (\mathbb{R}^n, g)$. Since the metric $g$ is uniformly equivalent to the Euclidean metric, we have $\mathcal{H}_d(L) = \mathcal{H}_d(g, \phi_g)$ for all $d \geq 0$. Our aim is to work with $h_d(g, \phi_g)$ instead. Let us first establish the following lemma.

LEMMA 1. *Let $\Gamma \subset \Omega$ be a smooth subdomain properly contained in an open domain $\Omega \subset M$. Let $\phi$ be a positive measurable function whose restriction to the boundary $\partial \Gamma$ of $\Gamma$ is also measurable. Suppose $H$ is a $k$-dimensional space of nonconstant $\phi$-harmonic functions defined on $\Omega$ and*

$$0 = \eta_1 \leq \eta_2 \leq \eta_3 \leq \ldots$$

*are the eigenvalues of $\partial \Gamma$ with respect to the operator $\phi^{-1} \operatorname{div}(\phi \bar{\nabla})$, where $\bar{\nabla}$ is the tangential gradient on $\partial \Gamma$. Then for any orthonormal basis $\{u_1, \ldots, u_k\}$ of $H$ with respect to the inner product*

$$D(u,v) = \int_\Gamma \langle \nabla u, \nabla v \rangle \phi \, dV,$$

*we have*

$$2 \sum_{i=1}^{k} \eta_i^{1/2} \leq \sum_{i=1}^{k} \int_{\partial \Gamma} |\nabla u_i|^2 \phi \, dA,$$

*where $dA$ is the area element on $\partial \Gamma$ induced by the metric $g$.*

*Proof.* Let us first point out that the right-hand side of the inequality is independent of the choice of the $D$-orthonormal basis. Let $\{w_1, \ldots, w_{k-1}\}$ be the first $k-1$ eigenfunctions on $\partial \Gamma$ corresponding to the eigenvalues $\{\eta_1, \ldots, \eta_{k-1}\}$. Using a similar argument as in Lemma 1.2 of [L-W2], we can find a $D$-orthonormal basis $\{u_1, \ldots, u_k\}$ of $H$ such that

$$\int_{\partial \Gamma} u_i w_j \phi \, dA = 0$$

for all $1 \leq j < i \leq k$. Hence by the variational principle, we obtain that

(1.5) $$\eta_i \int_{\partial \Gamma} u_i^2 \phi \, dA \leq \int_{\partial \Gamma} |\bar{\nabla} u_i|^2 \phi \, dA.$$

On the other hand, we have

$$2\sum_{i=1}^{k} \eta_i^{\frac{1}{2}} = 2\sum_{i=1}^{k} \eta_i^{\frac{1}{2}} \int_{\Gamma} |\nabla u_i|^2 \phi \, dV$$

$$= 2\sum_{i=1}^{k} \eta_i^{\frac{1}{2}} \int_{\partial \Gamma} u_i \frac{\partial u_i}{\partial \nu} \phi \, dA$$

$$\leq \sum_{i=1}^{k} \eta_i \int_{\partial \Gamma} u_i^2 \phi \, dA + \sum_{i=1}^{k} \int_{\partial \Gamma} \left(\frac{\partial u_i}{\partial \nu}\right)^2 \phi \, dA.$$

Applying (1.5) to the first term on the right-hand side and then combining the two terms, we conclude that

$$2\sum_{i=1}^{k} \eta_i^{\frac{1}{2}} \leq \sum_{i=1}^{k} \int_{\partial \Gamma} |\nabla u_i|^2 \phi \, dA. \qquad \square$$

## 2. $L$-harmonic functions

We are now ready to prove the main theorem. The manifold $M = (\mathbb{R}^n, g)$ is assumed to be the $\mathbb{R}^n$ endowed with the Riemannian metric satisfying (1.2) and (1.3). The function $\phi_g$ is defined by (1.4).

THEOREM 2. *Let $h_d = \dim \mathcal{H}_d(g, \phi_g)$ be the dimension of the space of polynomial growth $\phi_g$-harmonic functions of at most degree $d$. Then for any sequence of nonnegative numbers $0 = a_0 < a_1 < \cdots < a_j = d$, $h_d$ must satisfy*

$$\sum_{i=1}^{j} (a_i - a_{i-1}) h_{a_{i-1}} \leq \left(\frac{\Lambda_\infty}{\lambda_\infty}\right)^{n-1} \frac{2}{n!} (d + 2n - 1)^n.$$

*In particular, when $d$ is a positive integer, by taking $a_i = i$, we see that*

$$\sum_{i=1}^{d} h_i \leq \left(\frac{\Lambda_\infty}{\lambda_\infty}\right)^{n-1} \frac{2}{n!} (d + 2n)^n.$$

*Proof.* Let us define

$$\mathcal{H}'_d(g, \phi_g) = \{f \in \mathcal{H}_d(g, \phi_g) \mid f(0) = 0\}$$

and its dimension is denoted by $h'_d = \dim \mathcal{H}'_d(g, \phi_g)$. Obviously $h'_d = h_d - 1$, and we will estimate $h'_d$ instead. Let us denote the Dirichlet form on $B(r)$ by

$$D_r(f_1, f_2) = \int_{B(r)} \langle \nabla_g f_1, \nabla_g f_2 \rangle \phi_g \, dV_g.$$

Using the sequence $0 = a_0 < a_1 < a_2 \cdots < a_j = d$, we can decompose the space $\mathcal{H}'_d(g, \phi_g)$ with respect to the inner product $D_1$ into a direct sum

$$\mathcal{H}'_d(g, \phi_g) = H_1 \oplus \cdots \oplus H_j,$$

where each

$$H_i = \mathcal{H}_{a_i}(g, \phi_g) \ominus \mathcal{H}_{a_{i-1}}(g, \phi_g)$$

is the subspace consisting of harmonic functions of growth order between $a_{i-1}$ and $a_i$. Let us denote its dimension by $\dim H_i = k_i$. If we view the Dirichlet form $D_r$ as a bilinear form defined on the space $\mathcal{H}'_d(g, \phi_g)$, then the determinant of $D_r$ when computed with respect to a $D_1$-orthonormal basis must satisfy the growth assumption

(2.1) $$\det\nolimits_{D_1} D_r \leq C\, r^s,$$

where $s = \sum_{i=1}^{j} (2(a_i - 1) + n)\, k_i$.

For a fixed $r > 0$, and for any $\varepsilon > 0$, there exists a smooth matrix $(b^{ij})$ and a set $A \subset B(r)$ with its measure satisfying $|A| \leq \varepsilon$ such that

$$\lambda I \leq (b^{ij}) \leq \Lambda I$$

on $B(r)$, and

$$|a^{ij} - b^{ij}| < \varepsilon$$

on $B(r) \setminus A$. Moreover, for any fixed $r_0 < r$, we may also insist that

$$\lambda_{r_0} I \leq (b^{ij}) \leq \Lambda_{r_0} I$$

on $B(r) \setminus B(r_0)$.

Let us define the corresponding Riemannian metric $h = \sum_{i,j=1}^n h_{ij}\, dx_i\, dx_j$ on $\mathbb{R}^n$ by

$$(h^{ij}) = (h_{ij})^{-1} = w_h(x)\, (b^{ij})$$

with

$$w_h(x) = \sum_{i,j=1}^n b^{ij} \frac{\partial \rho}{\partial x_i} \frac{\partial \rho}{\partial x_j}.$$

We also define the corresponding function

$$\phi_h(x) = w_h^{\frac{n-2}{2}}(x)\, \sqrt{\det(b^{ij})}.$$

Note that we have

$$\lambda \leq w_h \leq \Lambda$$

and

$$\lambda^2 I \leq (h^{ij}) \leq \Lambda^2 I.$$

For each $u \in \mathcal{H}'_d(g, \phi_g)$, let $v \in H^1(B(r))$ be the solution to the Dirichlet problem

$$\sum_{i,j=1}^n \frac{\partial}{\partial x_i}\left(b^{ij}(x)\frac{\partial v}{\partial x_j}\right) = 0 \tag{2.2}$$

on $B(r)$ and

$$v = u \tag{2.3}$$

on $\partial B(r)$. We denote $\mathcal{K}(r)$ to be the space consisting of all such solutions $v$. The maximum principle then asserts that $\dim \mathcal{K}(r) = h'_d$. For $r_0 \leq t \leq r$, applying Lemma 1 to $H = \mathcal{K}(r)$ and $\Gamma = B(t)$, we get

$$2 \sum_{i=1}^{h'_d} \eta_i^{\frac{1}{2}}(\varepsilon, t) \leq \sum_{i=1}^{h'_d} \int_{\partial B(t)} |\nabla_h u_i|^2 \, \phi_h \, dA_h, \tag{2.4}$$

where $\eta_i(\varepsilon, t)$ are the eigenvalues of the operator $\phi_h^{-1} \mathrm{div}(\phi_h \bar{\nabla}_h)$ on $\partial B(t)$ with $\bar{\nabla}_h$ being the tangential gradient of the metric $h$ on $\partial B(t)$.

If we define the Dirichlet form on the space $\mathcal{K}(r)$ with respect to $h$ by

$$D_t(h)(f_1, f_2) = \int_{B(t)} \langle \nabla_h f_1, \nabla_h f_2 \rangle \, \phi_h \, dV_h,$$

then

$$(\ln \det\nolimits_{D_1(h)} D_t(h))' = \sum_{i=1}^{h'_d} \int_{\partial B(t)} \frac{|\nabla_h u_i|^2}{|\nabla_h \rho|} \, \phi_h \, dA_h.$$

Combining with (2.4) and the fact that

$$|\nabla_h \rho|^2 = h^{ij} \frac{\partial \rho}{\partial x_i} \frac{\partial \rho}{\partial x_j}$$
$$\leq \Lambda_{r_0}^2,$$

we conclude that

$$2 \sum_{i=1}^{h'_d} \eta_i^{\frac{1}{2}}(\varepsilon, t) \leq \Lambda_{r_0} \left(\ln \det\nolimits_{D_1(h)} D_t(h)\right)'.$$

Integrating from $r_0$ to $r$, we get

$$2 \int_{r_0}^r \sum_{i=1}^{h'_d} \eta_i^{\frac{1}{2}}(\varepsilon, t) \, dt \leq \Lambda_{r_0} \ln \det\nolimits_{D_{r_0}(h)} D_r(h). \tag{2.5}$$

We now claim that as $\varepsilon \to 0$, the solution $v \in H^1(B(r))$ to the boundary value problem (2.2) and (2.3) must satisfy

$$\int_{B(r)} |\nabla_g(u-v)|^2 \, dV_g \to 0. \tag{2.6}$$

Indeed,

$$\int_{B(r)} \sum_{i,j=1}^{n} b^{ij} \frac{\partial(u-v)}{\partial x_i} \frac{\partial(u-v)}{\partial x_j} dV_0$$

$$= -\int_{B(r)} (u-v) \sum_{i,j=1}^{n} \frac{\partial}{\partial x_i} \left( b^{ij} \frac{\partial(u-v)}{\partial x_j} \right) dV_0$$

$$= -\int_{B(r)} (u-v) \sum_{i,j=1}^{n} \frac{\partial}{\partial x_i} \left( b^{ij} \frac{\partial u}{\partial x_j} \right) dV_0$$

$$= -\int_{B(r)} (u-v) \sum_{i,j=1}^{n} \frac{\partial}{\partial x_i} \left( (b^{ij} - a^{ij}) \frac{\partial u}{\partial x_j} \right) dV_0$$

$$= \int_{B(r)} \sum_{i,j=1}^{n} \left( b^{ij} - a^{ij} \right) \frac{\partial(u-v)}{\partial x_i} \frac{\partial u}{\partial x_j} dV_0,$$

where $dV_0$ is the volume form with respect to the standard Euclidean metric. Splitting this integral over the two domains $B(r) \setminus A$ and $A$, we see that there exists a constant $\delta > 0$ with $\delta \to 0$ as $\varepsilon \to 0$, such that,

$$\int_{B(r) \setminus A} \sum_{i,j=1}^{n} \left( b^{ij} - a^{ij} \right) \frac{\partial(u-v)}{\partial x_i} \frac{\partial u}{\partial x_j} dV_0$$

$$\leq \varepsilon \int_{B(r)} |\nabla_0(u-v)| \, |\nabla_0 u| \, dV_0$$

$$\leq \varepsilon \left( \int_{B(r)} |\nabla_0(u-v)|^2 \, dV_0 \right)^{\frac{1}{2}} \left( \int_{B(r)} |\nabla_0 u|^2 \, dV_0 \right)^{\frac{1}{2}},$$

and

$$\int_{A} \sum_{i,j=1}^{n} \left( b^{ij} - a^{ij} \right) \frac{\partial(u-v)}{\partial x_i} \frac{\partial u}{\partial x_j} dV_0$$

$$\leq C \int_{A} |\nabla_0(u-v)| \, |\nabla_0 u| \, dV_0$$

$$\leq C \left( \int_{B(r)} |\nabla_0(u-v)|^2 \, dV_0 \right)^{\frac{1}{2}} \left( \int_{A} |\nabla_0 u|^2 \, dV_0 \right)^{\frac{1}{2}}$$

$$\leq C \delta \left( \int_{B(r)} |\nabla_0(u-v)|^2 \, dV_0 \right)^{\frac{1}{2}}.$$

Using the fact that the metrics $g$ and $h$ are uniformly equivalent to the Euclidean metric, we conclude that

$$\int_{B(r)} |\nabla_g(u-v)|^2 \, dV_g \leq C \int_{B(r)} |\nabla_h(u-v)|^2 \, dV_h$$

$$\leq C\, \delta \left( \int_{B(r)} |\nabla_0(u-v)|^2 \, dV_0 \right)^{\frac{1}{2}}$$

$$\leq C\, \delta \left( \int_{B(r)} |\nabla_g(u-v)|^2 \, dV_g \right)^{\frac{1}{2}},$$

and the claim is proved. In particular, we conclude from (2.6) that

$$(2.7) \qquad \int_{B(t)} |\nabla_h v|^2 \, dV_h \to \int_{B(t)} |\nabla_g u|^2 \, dV_g$$

as $\varepsilon \to 0$.

We now claim that for each $k \geq 1$, the $k^{\text{th}}$ eigenvalue $\eta_k(\varepsilon,t)$ must satisfy

$$(2.8) \qquad \eta_k(\varepsilon,t) \geq \lambda_{r_0}^2 \, \eta_k(1) \, t^{-2},$$

where $\eta_k(1)$ denotes the $k^{\text{th}}$ Euclidean eigenvalue of $\partial B(1)$ with respect to the induced Euclidean metric. In fact, the variational characterization implies

$$\eta_k(\varepsilon,t) = \inf_W \sup_{f \in W} \frac{\int_{\partial B(t)} |\bar\nabla_h f|^2 \, \phi_h \, dA_h}{\inf_{c \in \mathbb{R}} \int_{\partial B(t)} (f-c)^2 \, \phi_h \, dA_h},$$

where the infimum is taken over all $k$-dimensional subspace $W$ in the space of smooth functions on $\partial B(t)$. By the co-area formula, we have

$$|\nabla_h \rho|^{-1} \, dA_h = \sqrt{\det(h_{ij})} \, dA_0,$$

where $dA_h$ and $dA_0$ are the area elements on the sphere induced by the metric $h$ and the Euclidean metric respectively. Thus, by the choice of $\phi_h$, we have

$$dA_0 = \phi_h \, dA_h.$$

So,

$$\eta_k(\varepsilon,t) = \inf_W \sup_{f \in W} \frac{\int_{\partial B(t)} |\bar\nabla_h f|^2 \, dA_0}{\inf_{c \in \mathbb{R}} \int_{\partial B(t)} (f-c)^2 \, dA_0}$$

$$\geq \lambda_{r_0}^2 \, \inf_W \sup_{f \in W} \frac{\int_{\partial B(t)} |\bar\nabla_0 f|^2 \, dA_0}{\inf_{c \in \mathbb{R}} \int_{\partial B(t)} (f-c)^2 \, dA_0}$$

$$= \lambda_{r_0}^2 \, \eta_k(t)$$

$$= \lambda_{r_0}^2 \, \eta_k(1) \, t^{-2},$$

where $\eta_k(t)$ is the $k^{\text{th}}$ Euclidean eigenvalue of the sphere $\partial B(t)$. This proves the assertion (2.8).

In is known that [B-G-M]) for each $k \geq 1$,

$$\eta_k(1) = q^2 + (n-2)q$$

where $q$ is the least positive integer satisfying

(2.9) $$k \leq \binom{n+q-1}{q} + \binom{n+q-2}{q-1}.$$

However, one checks easily that

$$k \leq \frac{2}{(n-1)!}\left(q + \frac{n-1}{2}\right)^{n-1};$$

hence we have

$$q \geq \left(\frac{(n-1)!}{2}k\right)^{\frac{1}{n-1}} - \frac{n-1}{2},$$

and

$$\eta_k^{\frac{1}{2}}(1) \geq \left(\frac{(n-1)!}{2}k\right)^{\frac{1}{n-1}} - \frac{n-1}{2}.$$

Therefore,

(2.10) $$\sum_{i=1}^{k} \eta_i^{\frac{1}{2}}(1) \geq \left(\frac{(n-1)!}{2}\right)^{\frac{1}{n-1}} \sum_{i=1}^{k} i^{\frac{1}{n-1}} - \frac{n-1}{2}k$$

$$\geq \left(\frac{(n-1)!}{2}\right)^{\frac{1}{n-1}} \int_0^{k-1} x^{\frac{1}{n-1}}\, dx - \frac{n-1}{2}k$$

$$= \left(\frac{(n-1)!}{2}\right)^{\frac{1}{n-1}} \frac{n-1}{n}(k-1)^{\frac{n}{n-1}} - \frac{n-1}{2}k$$

$$\geq \left(\frac{(n-1)!}{2}\right)^{\frac{1}{n-1}} \frac{n-1}{n}k^{\frac{n}{n-1}} - (n-1)k.$$

Applying (2.8), (2.10) to (2.5), we have

(2.11) $$2\lambda_{r_0}\left(\left(\frac{(n-1)!}{2}\right)^{\frac{1}{n-1}} \frac{n-1}{n}(h_d')^{\frac{n}{n-1}} - (n-1)h_d'\right)(\ln r - \ln r_0)$$

$$\leq 2\int_{r_0}^{r} \sum_{i=1}^{h_d'} \eta_i^{\frac{1}{2}}(\varepsilon, t)\, dt$$

$$\leq \Lambda_{r_0} \ln \det{}_{D_{r_0}(h)} D_r(h).$$

Letting $\varepsilon \to 0$ and using (2.7), we conclude from (2.11) that

$$2\left(\frac{\lambda_{r_0}}{\Lambda_{r_0}}\right)\left(\left(\frac{(n-1)!}{2}\right)^{\frac{1}{n-1}} \frac{n-1}{n}(h'_d)^{\frac{n}{n-1}} - (n-1)h'_d\right)(\ln r - \ln r_0)$$
$$\leq \ln \det_{D_{r_0}(g)} D_r(g)$$
$$\leq \left(\sum_{i=1}^{j}(2(a_i-1)+n)k_i\right)(\ln r + C),$$

where we have used (2.1) and the constant $C$ is independent of $r$. Dividing both sides by $\ln r$ and letting first $r \to \infty$ and then $r_0 \to \infty$, we get

$$2\left(\frac{\lambda_\infty}{\Lambda_\infty}\right)\left(\left(\frac{(n-1)!}{2}\right)^{\frac{1}{n-1}} \frac{n-1}{n}(h'_d)^{\frac{n}{n-1}} - (n-1)h'_d\right)$$
$$\leq \sum_{i=1}^{j}(2(a_i-1)+n)k_i.$$

Using the fact that $k_i = h'_{a_i} - h'_{a_{i-1}}$, we can rewrite this inequality into the form

$$(2.12) \quad \sum_{i=1}^{j}(a_i - a_{i-1})h'_{a_{i-1}} \leq \left(d - 1 + \frac{n}{2}\right)h'_d$$
$$- \left(\frac{\lambda_\infty}{\Lambda_\infty}\right)\left(\left(\frac{(n-1)!}{2}\right)^{\frac{1}{n-1}} \frac{n-1}{n}(h'_d)^{\frac{n}{n-1}} - (n-1)h'_d\right)$$
$$\leq \left(d + \frac{3n}{2} - 2\right)h'_d - \left(\frac{\lambda_\infty}{\Lambda_\infty}\right)\left(\frac{(n-1)!}{2}\right)^{\frac{1}{n-1}} \frac{n-1}{n}(h'_d)^{\frac{n}{n-1}}.$$

As a function of $h'_d$, the right-hand side of (2.12) achieves its maximum when $h'_d$ satisfies

$$\left(\frac{\lambda_\infty}{\Lambda_\infty}\right)\left(\frac{(n-1)!}{2}\right)^{\frac{1}{n-1}}(h'_d)^{\frac{1}{n-1}} - \left(d + \frac{3n}{2} - 2\right) = 0.$$

Plugging this value of $h'_d$ into (2.12) and simplifying, we obtain

$$\sum_{i=1}^{j}(a_i - a_{i-1})h'_{a_{i-1}} \leq \frac{2}{n!}\left(\frac{\Lambda_\infty}{\lambda_\infty}\right)^{n-1}\left(d + \frac{3n}{2} - 2\right)^n.$$

Note that $h'_d = h_d - 1$; thus we have

$$\sum_{i=1}^{j}(a_i - a_{i-1})\, h_{a_{i-1}} = \sum_{i=1}^{j}(a_i - a_{i-1})\, h'_{a_{i-1}} + d$$

$$\leq \frac{2}{n!}\left(\frac{\Lambda_\infty}{\lambda_\infty}\right)^{n-1}\left(d + \frac{3n}{2} - 2\right)^n + d$$

$$\leq \frac{2}{n!}\left(\frac{\Lambda_\infty}{\lambda_\infty}\right)^{n-1}(d + 2n - 1)^n.$$

This proved the first assertion. In the case $d$ is a positive integer and $a_i = i$, the first assertion implies that

$$\sum_{i=1}^{d} h_{i-1} \leq \frac{2}{n!}\left(\frac{\Lambda_\infty}{\lambda_\infty}\right)^{n-1}(d + 2n - 1)^n.$$

Now replacing $d$ by $d+1$, we get the second assertion. $\square$

COROLLARY 3. *Let*

$$Lf = \sum_{i,j=1}^{n} \frac{\partial}{\partial x_i}\left(a^{ij}(x)\frac{\partial f}{\partial x_j}\right)$$

*be a uniformly elliptic operator defined on $\mathbb{R}^n$ with measurable coefficients satisfying* (0.1). *If $h_d(L) = \dim \mathcal{H}_d(L)$ denotes the dimension of the space of polynomial growth L-harmonic functions of at most degree $d$. Then*

$$\sum_{i=1}^{d} h_i(L) \leq \left(\frac{\Lambda_\infty}{\lambda_\infty}\right)^{n-1}\frac{2}{n!}\left(d^n + O(d^{n-1})\right)$$

*and*

$$\liminf_{d\to\infty} d^{-(n-1)}\, h_d(L) \leq \left(\frac{\Lambda_\infty}{\lambda_\infty}\right)^{n-1}\frac{2}{(n-1)!}.$$

UNIVERSITY OF CALIFORNIA, IRVINE, CA
*E-mail address*: pli@math.uci.edu

UNIVERSITY OF MINNESOTA, MINNEAPOLIS, MN
*E-mail address*: jiaping@math.umn.edu